\documentclass[preprint,number,8pt]{elsarticle}
\usepackage{nonfloat}
\usepackage{amsmath}
\usepackage{amsfonts}
\usepackage{amssymb}
\usepackage{graphicx}
\newcommand{\s}{$A_s$}
\newcommand{\ns}{$A_{ns} \,$}
\newcommand{\nd}{$A_{nd}\,$}
\newcommand{\bi}{$B_{int}\,$}
\newcommand{\abs}{$B_{abs}\,$}
\usepackage{ucs}
\usepackage{times}
\usepackage{epsfig}
\usepackage{epstopdf}
\usepackage[active]{srcltx}
\usepackage[mathscr]{eucal}
\usepackage{caption}
\usepackage{url}
\usepackage{color}
\journal{Theoretical biology}
\begin{document}
\begin{frontmatter}
\title{Mathematical Modeling of Transport and Degradation of Feedstuffs in the Small Intestine}
 \author{Masoomeh Taghipoor$^{\dagger *}$, Philippe Lescoat$^{\dagger}$, Jean-Ren\'e Licois$^*$, Christine Georgelin$^*$, Guy Barles$^*$}
\address{$*$ Laboratoire de Math\'ematiques et Physique Th\'eorique\\
Unit\'e Mixte de Recherche $6083$ du CNRS\\
F\'ed\'eration de Recherche $2964$ Denis Poisson\\
Universit\'e Fran\c{c}ois Rabelais, Parc de Grandmount\\
37200 Tours, France\\
$\dagger$ INRA, UR83 Recherches Avicoles, 37380 Nouzilly, France }

\begin{abstract}
We describe a mathematical model of digestion in the small  
intestine. The main interest of our work is to consider simultaneously the different aspects of digestion i.e. transport of  
the bolus all along the intestine, feedstuffs degradation according to  
the enzymes and local physical conditions, and nutrients absorption. A  
system of coupled ordinary differential equations is used to model  
these phenomena. The major unknowns of this system are the position of  
the bolus and its composition. This system of equations is solved  
numerically. We present several numerical computations for the  
degradation, absorption and transport of the bolus with acceptable  
accuracy regarding the overall behavior of the model and also when challenged versus experimental data. The main feature and interest of  
this model are its genericity. Even if we are at an early stage of  
development, our approach can be adapted to deal with contrasted feedstuffs in non-ruminant  
animal to predict the composition and velocity of bolus in the small  
intestine.
\end{abstract}
\end{frontmatter}

\section{Introduction}\label{introduction}
The main step of digestion and absorption along the gastrointestinal tract takes place in the small intestine for non ruminant animals. To reach an optimized composition of available nutrients due to their behavior in the digestive system, it is necessary to understand and predict the digestion and absorption of the ingested feedstuffs in the small intestine \cite{karasov,williams,logan2}. 
It is also now well-known that the use of implanted experimental devices may modify the dynamic of digestion in the small intestine \cite{Srivastava, wilfart}.

Several models have been developed representing the digestion and transport of bolus in the small intestine. In the model of \cite{bestianelli} digestion and absorption are integrated and represented considering only the polymers and individual absorbable end products. The transit through the small intestine is modeled mainly as a result of gastric emptying. No peristaltic wave is taken into account, and the bolus contained only the dry matter. \cite{logan} describe the digestion and absorption using the plug flow reactors to encapsulate complex digestion phenomena in a simple set of equations. Different rate of absorption and degradation are involved: first order kinetics, Michaelis-Menten and the sigmoid ones. A detailed model of the intestinal propulsion is provided by \cite{miftahof,Pullana}. However, these models portray the transport of bolus simplistically, or they represent only a limited number of different processes involved in digestion.

This article tries to go further in the modeling of digestion in the small intestine by considering the different steps of digestion i.e. the transport of the bolus all along the intestine, feedstuffs degradation according to the enzymes and local physical conditions, and nutrient absorption. Therefore a system of coupled ordinary differential equations is used. The major unknowns of this system are the position of the bolus and its composition.

In fact, several models are presented reflecting the modeling process at its different stages with our attempts to make it more realistic by inclusion of more sophisticated and relevant biological phenomena and chemical transformations. We decided to describe the different steps with the assumptions leading us to our choices instead of presenting only the last model since the whole process by itself may help to underline relevant questions to be further discussed. Of course, this modeling process is an iterative one and is still going ahead in directions which are described in Section~\ref{CandP}.

Our models intend to be a mechanistic approach of feedstuffs digestion even though simplifications were included according to participatory approaches between biologists and mathematicians. However they involve a lot of different unknowns and parameters, and require a numerical software to obtain suitable approximation of the solutions since it is hopeless to obtain explicit ones. Scilab software was used to perform these numerical computations \footnote{The reader can perform its own numerical experiments, with the possibility of changing
the parameters, using our Scilab software online at the URL :
\url{http://www.lmpt.univ-tours.fr/modingre}}
. 

In all our models, we try to estimate the parameters using data from scientific literature. When these data are not available, we assume the reasonable values for the parameters.

The article is organized as follows : Section~2 is devoted to present the main assumptions of our models and most of our notations. In Section~3, we describe the transport equations. In our four different models, we point out that there are only two different ways of modeling the transport of the bolus in the intestinal tract. The main differences concern the degradation itself, with several possibilities for the composition of the bolus, for the enzymatic reactions and the water influence. The outcoming stages (4 different models) are presented in Section~4, with the key assumptions and characteristics of each model. Section~5 is a comparison of these four models and of the numerical results of the most sophisticated model (Model~3) versus some experimental data from the literature.  Finally, in Section~6, we criticize our models and describe the perspective.

\section{General Hypothesis and Synthetic Presentation of the Different Models} 

\begin{figure}[ht]  
\begin{center}  
\includegraphics[width=90mm]{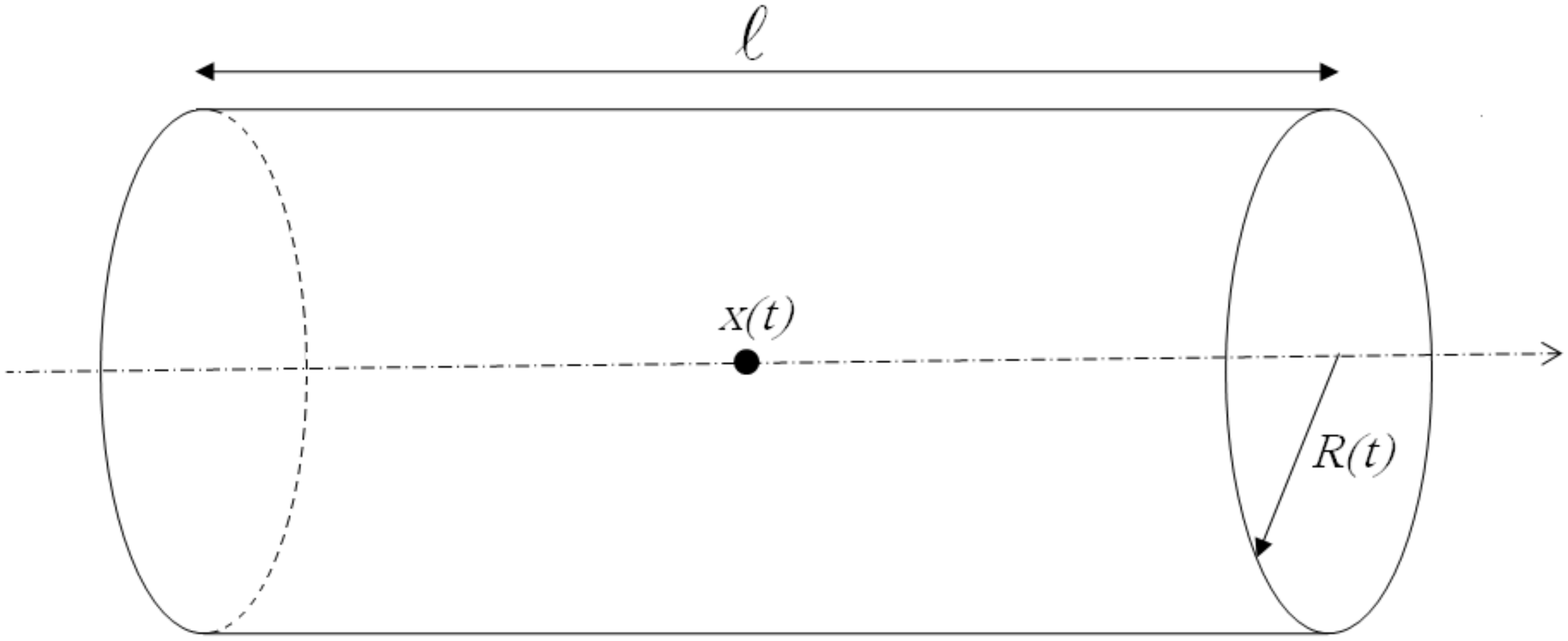}
\end{center} 
\caption{  \label{fig1.png}
The Cylindric bolus with its different characteristics} 
\end{figure}

Common assumptions to all models are the following ones. 
\begin{itemize}
\item[$(i)$] The first simplification concerns the small intestine representation itself. Instead of taking into account its complex geometry, it is represented as a one-dimensional interval $[0,L]$. The position of the bolus in the small intestine at time $t$ is given by $x(t) \in [0,L]$ (cf. figure\ref{fig1.png}). The equation of the transport of bolus along the small intestine is described on $x(t)$. 

\item[$(ii)$] The bolus is treated as a homogeneous cylinder with a fixed length $\ell $ and variable radius $R(t)$ which is a function of time $t$. To locate this cylinder, we use the position $x(t)$ of its center. This assumption is justified by the general shape of the small intestine's segments as well as the observation of the real bolus in animals' small intestines. As the length of bolus is assumed fixed, the degradation only changes the radius. Degradation of substrates is obtained by enzymatic reactions with possible subsequent absorption by the intestinal wall \cite{digestion}.
\item[$(iii)$] The enzymes which participate in enzymatic hydrolysis included in our models consist systematically in pancreatic and brush-border ones with the possibility of including exogenous and gastric ones. The enzymes' activity depends on the \emph{pH} of the small intestine at each point along its length. The brush-border~enzymes~on the intestinal wall are assumed to be always in excess \cite{digestion}. 
\end{itemize}	

The bolus moves through the intestinal tract because of the pulses resulting from peristaltic waves and gastric emptying which gives an initial velocity to the bolus \cite{digestion}.
Peristalsis are series of wavelike contractions occurring in the smooth muscle layer of the gastrointestinal tract. It is a physiological process that results in intestinal motility and propulsion of ingested food along the intestine. It starts as a ring-like constriction initially which later moves mostly forward along the intestine. 
Moreover it might be assumed that it helps the bolus to be digested by spreading the food particles along the intestinal wall for effective digestion and absorption.
Therefore to model the transport, bolus movement is connected to pulses all along the intestine with the initial velocity coming from the gastric emptying effect. Efficiency of pulses is proportional to the volume of the bolus and it is inversely correlated to the distance between the bolus and the pylorus \cite{J.Fioramonti,Xiao}. The bolus speed is assumed to be slowed down by the local conditions in the small intestine lumen (friction on the borders, viscosity effects,...etc). The effects of these different local factors depends on the composition of the bolus, and in particular its dilution.
 
For the bolus content, the following assumptions and notations are used throughout this paper
\begin{itemize}
\item[$\bullet$] The bolus includes a single species whose total mass is denoted by $A$. In the most sophisticated model presented in this work, $A$ is composed of $A_s$, $A_{ns}$ and $ A_{nd} $, in other words $A=$\s +\ns + \nd. $A_s$ is the mass of solubilized fraction of $A$ which can be hydrolyzed in the presence of the enzymes. The index $s$ stands for ``solubilized``. $A_{ns}$ is the mass of ``non-solubilized`` fraction of $A$, for example dry starch or the non-emulsified lipids. Transformation of $A_{ns}$ into $A_s$ requires a sufficient quantity of water. Regarding lipids emulsion, we assume that the bile salts are in excess. The mechanism $A_{ns} \leftrightarrow A_s$ is described through an equilibrium property depending on the water quantity in the bolus. Finally $A_{nd}$ is the mass of non-degradable $A$, which enters and leaves the small intestine without any change. For example the vegetal fiber in feed are poorly digestible. Moreover the fiber matrix of feedstuffs or the anti-nutritional factor content can be responsible for a reduction in the digestibility of some amino acids in some feedstuffs \cite{bernard}.

\item[$\bullet$] The quantity $B$ is the mass of product obtained from $A_s$ by enzymatic reactions, it is composed of $B_{int}$ and $B_{abs}$
$$ A_s\big(\stackrel{enzymes}{\longrightarrow}B_{int}\big)\stackrel{enzymes}{\longrightarrow} B_{abs}.$$
The quantity $B_{int}$ is the product of hydrolysis due to gastric and pancreatic enzymes, the index $int $ stands for ''intermediate'' substrate which is not yet absorbable since it is not fully degraded. This transformation has to be completed by a second one at the border of the small intestine via the brush-border enzymes (e.g. : proteins being degraded to polypeptides and afterwards to dipeptide or amino acids, which are absorbable). This second transformation is also able to give \abs  directly from $A$. The quantity $B_{abs}$ is the absorbable fraction with index $abs$ indicating ``absorbable``.
\item[$\bullet$]The quantity $e$ represents the gastric and pancreatic enzymes.
\item[$\bullet$]The quantity $W$ is the mass of water in the bolus and $[W]$ indicates the proportion of water in the bolus : $W/(A+B+W)$.
\item The quantity $V(t)$ denotes the volume of bolus which is equal to $(A+B+W)/\rho $, where $\rho $ denotes the density. For the sake of simplicity, we assume that all the substrates of bolus have the same density $\rho$. The total mass of the bolus is $(A+B+W)(t)$ at each moment.
\end{itemize}
Digestion consists in the transformation of digesta to absorbable nutrients through enzymatic hydrolysis. Volumic transformation is the degradation of $A_s$ into $B_{int}$ inside the bolus and transformation on the bolus surface is the degradation of both $A_s$ and $B_{int}$ into $B_{abs}$ on a thin layer around the bolus.
The following hypothesis are added progressively with upgraded versions of the model.

\begin{itemize}
\item[Model 1.] In the first model, the whole bolus is considered to be solubilized ($A=$\s). $A$ is hydrolyzed thanks to gastric enzymes and becomes directly absorbable nutrients ($B=$\abs). In this model, brush-border or pancreatic ones are not taken into account. Such mechanisms are associated for example with the consumption of disaccharides (resp. monosaccharides) such as milk sugar(resp. glucose).
%$$ A_s \rightarrow B_{abs} .$$
\item[Model 2.] This model is an attempt to have a more realistic modeling of degradation. The bolus is still assumed to be completely solubilized. The absorbable nutrients can be obtained by two ways: either by a direct surfacic transformation $A \rightarrow B_{abs}$ or through a first volumic degradation $A\rightarrow B_{in}$ followed by a second one $B_{int} \rightarrow  B_{abs}$ at the bolus surface by brush-border enzymes.
\item[Model 3.] This model includes the solubilization of the bolus in presence of water. A is splitted into $A_s$, \ns, \nd. Equations are added to express the equilibrium $A_s \leftrightarrow A_{ns}$ which depends on the quantity of water. The non-degradable part of bolus enters and leaves the small intestine without any mechanical or chemical change in its initial form. A key feature of this model concerns the transport of bolus along the small intestine since it is connected to the quantity of water through lubrification effects.	
\item[Model 4.] This model is a simplification of the previous one by mathematical arguments. Through homogenization methods it is shown that the acceleration can be averaged and an equation with this averaged acceleration is substituted for the pulses in transport equation. 
Detailed models are described in Section $4$.
\end{itemize}
\section{Transport}\label{transport}
We present a mathematical formulation of the transport of bolus in the small intestine. It is based on the physiology of the pig's small intestine to get consistent parameters.
The duodenum is characterized by oscillatory electrical events (slow waves) occurring at a rate of 18/min. After food intake, some of these waves are associated with spikes bursts which are responsible for contractions and therefore propelling the bolus through the small intestine. We assume only 6 of the 18 slow waves by minute are followed by the spikes. It leads to one efficient contraction every 10 seconds \cite{J.Fioramonti}.
 %After food intake, frequency of the peristaltic waves initiated in duodenum is about 18 per minute. 
The mean transit time of each peristaltic wave is assumed to be 150 minutes to move along small intestine from duodenum to the end of ileum according to \cite{Laplace} and \cite{Rayner}. We assume also that the pig's small intestine is about 18 meters \cite{Pommier}, hence the average velocity $c$ of these waves is $7,2$ m/h.
Each peristaltic wave takes $x(t)/c$ seconds to reach the bolus in position $x(t)$, therefore the pulse which pushes the bolus in time $t$ is generated in duodenum at time $t -x(t)/c.$

If $v(t)$ denotes the velocity of the bolus ($v(t)=\displaystyle \frac{dx}{dt}(t)$), the effect of pulses is modeled through the following equation   
\begin{equation}\label{equ:transport1} \nonumber
 \frac{d^2}{dt^2}x(t) = \frac{d}{dt}v(t)=\frac{d}{dt} y\left(t-x(t)/c\right),
\end{equation}
 the term  $\displaystyle \frac{d}{dt}y$  represents the pulses which are defined as a periodic function of period $10$ seconds such that $$\int^{10}_{0} \frac{d}{dt} y(t)dt=1$$ and for $t<0$, we assume $y^{'}(t)=0$.
 
Over a period, each pulse is an approximation of a Dirac mass of the origin. Therefore we define it as a function with the value $1/\epsilon$ during a very short interval of time $\epsilon $ and $0$ at all other time.

According to \cite{Xiao} and \cite{Rivest} the efficiency of the~peristaltic~waves~increases with the size of the bolus and decreases with the distance from pylorus. We assume that all these dependences are affine, namely 
\begin{equation}\label{equ:tansport}\nonumber
 \dfrac{d^2x}{dt^2}(t) =\dfrac{d}{dt}\left[ y\left(t- x(t)/c\right)\right]\dfrac{c_0 + c_1 V(t)}{a+bx(t)},
 \end{equation}
where, $c_0$ and  $c_1$ are determined under the assumption that the acceleration depends linearly on $V(t)$. THe constants $a$ and $b$ are obtained from experimental data. 

The intestinal lumen is a confined environment which prevents the bolus to move perfectly according to the previous equation : the bolus has to work its way through the small intestine and is also submitted to the friction with the intestinal wall. All these friction effects are related to the ``viscosity'' of the bolus and we have two different ways to model the friction term : either as a constant effect which is independent of the bolus composition (models $1$ \& $2$) or with a lubrification effect coming from the proportion of water in the bolus (models $3$ $\&$ $4$). More specifically, in Equation~(\ref{equ:transport2}) below, the coefficient $K(t)$ is either constant in models $1$ \& $2$ or $\tilde K/[W]$ in models $3$ $\&$ $4$, where $\tilde K$ is a constant.

The final equation of transport is therefore the following one 
\begin{eqnarray}\label{equ:transport2}
&\dfrac{d^2x}{dt^2}(t) &= \frac{d}{dt}[y(t-x(t)/c)]\frac{c_0+c_1V(t)}{a+bx(t)}-K(t) \frac{dx}{dt}(t)\\
&\dfrac{dx}{dt}(0) &= v_0 ,\, \, \,  x(0)=0 \nonumber 
\end{eqnarray}

\section{Digestion}\label{digestion}
Digestion is a mechanical and chemical process by which the feedstuffs molecules are broken down to the smaller ones by enzymes in order to to make them available for absorption. The uptake of the obtained nutrients is mainly by absorption.
In this section the different steps of modeling are detailed.
 \subsection{Model 1}
%%%hypotheses%%%%%%%%%%%%%
In this model the bolus is assumed to be completely solubilized ($A$=\s). We also assume that the necessary enzymes for hydrolysis are mixed with the bolus in the stomach. The product of following reaction is directly absorbable ($B=$\abs)
$$A+e \rightarrow B.$$
$e$ denotes the gastric enzymes.
%%%%%%%%%%Inconnues%%%%%%%%%%%
The first aim of this model is to define the variation of the bolus which means the amount of the different substrates $A$, $B$ and $e$ at every time in the cylinder, and the second is to locate the bolus along the small intestine.
%%%%equation digestion%%%
We assume that, the evolution of $A$ or its volumic transformation depends on its mass at each moment and the enzyme activity. This equation follows the law of mass action
\begin{equation} \nonumber
\frac{dA}{dt} =-Ck(x,e)A
\end{equation}
 where $C$ denotes the degradation rate and, $k(x,e)$ is the enzyme activity which depends on the $pH$ of the small intestine and the presence of the enzymes at each point along it.

 The product $B$ of the volumic transformation of $A$ is absorbed by intestinal wall with a constant rate $k_{abs}$
\begin{equation} \nonumber
\frac{dB}{dt}=Ck(x,e)A-k_{abs}B. \label{eq2}
\end{equation}
%%%%enzymes%%%%%%%%%
There are also the degradation and inactivation of the enzymes along the small intestine
\begin{equation}\nonumber
 \frac{de}{dt}=-k_e e 
\end{equation}
where $k_e$ is the rate of degradation of the enzymes which depends on their types.  
The activity of each enzyme as a function of $pH$ of small intestine is roughly known. We know also the $pH$ of each point along the small intestine. The composition of these two functions gives the enzyme activity at each point $x$ along it. 
%%%%%%%%%%%%%%%%%%%%%%%%%%%%%%%%%%%%%%%%%%%%%%%%%%%%%%%%%%%%%%%%%%%%%%%%
\subsection{Model 2}

In this second model, the presence of pancreatic enzymes in the small intestine as well as the brush-border ones on its wall are considered.
 The pancreatic secretions help neutralizing the stomach acid as they enter the small intestine. They also contain pancreatic enzymes. The level of the secretions is a function of volume and composition of the bolus entered the small intestine. The brush border enzymes are the enzymes for the terminal stage of digestion which is the surfacic hydrolysis. Contrary to the pancreatic enzymes they are not free in the intestinal lumen, but rather, in the plasma membrane of the enterocyte. 

%%%%Hypotheses%%%%%%%
We assume that the bolus is completely solubilized. The product of the hydrolysis $B$ consists in \bi and \abs ($B=$\bi +\abs ).

The following scheme represents the chemical reactions of the bolus in this model
$$ A_s\rightarrow B_{int},\,\,\quad A_s  \rightarrow   B_{abs}, \,\,\quad  B_{int}  \rightarrow  B_{abs} .$$
The first reaction takes place inside the bolus by pancreatic and gastric enzymes, the second and the third ones take place on the surface of the bolus.
%%%%inconnus%%%

%%%%%%%%%%degradation%%%%%%%%%%%%%%%%%%%%%%%%%%%%%%%%%%%%%%%%%%%%%%%%%%%%%%%%%%%%%%%%%%%%%%%%%%%%%%%
The degradation of $A$ in this model is the result of the volumic hydrolysis of $A$ as in Model~1, and its surfacic hydrolysis by brush-border enzymes 

\begin{eqnarray} \nonumber
\frac{dA}{dt}&=&-Ck(x,e)A-C_{abs} (2 \pi R \ell)\frac{A}{A+B_{int}+B_{abs}}\\\nonumber \label{eq2d}
& = &-Ck(x,e)A -2C_{abs}\sqrt{\pi l/\rho} \frac{A}{(A+B_{int}+B_{abs})^{1/2}}, \nonumber
\end{eqnarray}
the second term represents surfacic transformation of $A$ to \abs. We recall that the mass of the bolus in this model is
$$A(t)+B_{int}(t)+B_{abs}(t)=\rho V(t)=\rho \pi R^2(t) l ,$$
and therefore the lateral surface of the cylinder is given by
$$ 2 \pi R \ell = 2 \sqrt{\pi l/\rho} (A+B_{int}+B_{abs})^{1/2}\; .$$

This transformation depends on the fraction of $A$ on the surface of the bolus which is written by $(2 \pi R \ell)\dfrac{A}{A+B_{int}+B_{abs}}$. The unit of the degradation coefficient per unit of surface and time, $C_{abs}$ , is $g.m^{-2}.s^{-1} $.

After a distance traveled  by bolus of about $5\%$ of the total length of the small intestine which is approximatively $85\, cm$ in an growing pigs, the input of secretions starts and it stops after a distance of $\alpha$ meters traveled by bolus. We assume their mass is about $\beta\%$  of the bolus mass. In the following equation, the effect of these secretions on the variation of $A$ is taken into account 
\begin{eqnarray*}
\frac{dA}{dt} &=&\text{...}+\ln(1.\beta)\frac{1}{\alpha} \frac{dx}{dt} \chi\left(\left(x(s)-0.85\right)/\alpha\nonumber \right)A,
\end{eqnarray*}
where $\chi $ is a localization function in the above equation which reflects the fact that secretions arrive in the small segment of the intestine, say between $0.85$ cm and $0.85 + \alpha$ cm.\\

%%%%%%bint%%%%%%%%%%%%%%%%%%%%%%%%%%%%%%%%%%%%%%%%%%%%%%%%%%%%%%%%%%%%%%%%%%%%%%%%%%%%%%%%%%%%%
The product of volumic hydrolysis, \bi , participates in the creation of \abs on the surface of the bolus. Therefore its variation is modeled by

\begin{eqnarray*}
\frac{dB_{int}}{dt} &=&Ck(x,e)A+\ln(1.25)\frac{1}{\alpha}\frac{dx}{dt} \chi ((x(s)-0.85)/\alpha)B_{int} \\ 
&& -2{ C_{iabs}}\sqrt{\pi l/ \rho }\frac{B_{int}}{(A+B_{int}+B_{abs})^{1/2}} .\label{eq:bi}
\end{eqnarray*}
%%%%%%%%%absorbable%%%%%%%%%%%%%%%
The absorbable nutrients on the bolus are not absorbed instantaneously \cite{logan}. In this model we assume that the absorption rate follows Michaelis-Menten mechanism. The constant $k_{abs} $ is the maximal rate of absorption at saturation, $k$ is the Michaelis constant which is half saturation 
\begin{equation*} 
\frac{dB_{abs}}{dt}=2\sqrt{\pi l/ \rho }\frac{C_{abs}A+{C_{iabs}}B_{int}}{(A+B_{int}+B_{abs})^{1/2}}-k_{abs}\frac{B_{abs}}{k+B_{abs}}.\\ \label{eq2a}
\end{equation*}

%Model3 %%%%%%%%%%%%%%%%%%%%%%%%%%%%%%%%%%%%%%%%%%%%%%%%%%%%%%%%%%%%%%%%%%%%%%%%%%%%%%%%%%%%%%

\subsection{Model 3}

In this model the ingested food consists in $A_{ns}$, $A_{nd}$, $A_s$ and water ($A=A_{ns}+ A_{nd}+ A_s+W$). We incorporate two effects of water on digestion~: the first one is the dilution of the bolus and its impacts on degradation and absorption and the second one is the lubrification and its consequences on the transport.
  
We assume that the evolution of $A_s$ and $A_{ns}$ aims at reaching an equilibrium in which the ratio between $A_s$ and $A_{ns}$ is fixed and depends only on the proportion of water, namely $A_{s}~=~ \mu \left([W]\right) A_{ns} $ stressing that solubilization of $A_{ns}$ depends on bolus dilution. From the mathematical standpoint, we write this evolution as
\begin{equation}\label{eq:ns}
\frac{dA_{ns}}{dt}=-k_s \bigg(\mu \left([W]\right) A_{ns} - A_{s}\bigg), 
\end{equation}
where $\mu $ is a linear function of water and the constant $k_s$ represents the return rate to equilibrium.

The amount of water in the intestinal lumen is regulated by several complex biological phenomena. 
In fact the proportion of water in the bolus aims at reaching $[W_0]$ in a rather fast way which we translate it on a mathematical standpoint
 \begin{equation}
 \frac{d[W]}{dt}=-k_w([W]-[W_0])+\ln(1.\beta) \frac{1}{\alpha}\frac{dx}{dt}\chi ((x(s)-0.85)/\alpha) [W]
\end{equation}
where $k_w$ is large enough to reach the equilibrium in an adequate time. The second term of above equation is the fraction of water in pancreatic secretions.

The variation of $A_s$ depends on its degradation by volumic and surfacic hydrolysis, and contribution of pancreatic secretions  as in previous models. It also depends on the equilibrium with \ns resulting from equation (\ref{eq:ns}) which is the first term of equation below
\begin{equation}\label{eq3d}
\begin{split}
\frac{dA_s}{dt}&=k_s \bigg(\mu \left([W]\right) A_{ns} - A_{s}\bigg)-Ck(x,e)A_s(t)\\
&-2C_{abs}\sqrt{\pi l/\rho} \frac{A_s}{(A_s+A_{ns}+A_{nd}+B_{int}+W+B_{abs})^{1/2}}[W] \\
&+\ln(1.25)\frac{1}{\alpha}\frac{dx}{dt}\chi\left((x(s)-0.85\right)/\alpha)A_s.
\end{split}
\end{equation}
\\
The variation of absorbable nutrients depends on the creation of \abs by enzymatic hydrolysis of  \s  and  \bi and its absorption by intestinal wall

\begin{align}
\frac{dB_{abs}}{dt}= & 2\sqrt{\pi l/ \rho }\frac{C_{abs}A+{C_{iabs}}B_{int}}{(A_s+A_{ns}+A_{nd}+B_{int}+W+B_{abs})^{1/2}}[W] \nonumber \\
& \ \qquad -k_{abs}\frac{B_{abs}}{k+B_{abs}}.\label{eq3a}
\end{align}
The non-degradable fraction of $A$, namely $A_{nd}$, enters in the small intestine and leaves it without any change in its structure.

As we already indicated in Section~\ref{transport}, lubrification of the bolus~depends~on the presence of water. For this model, the friction coefficient in equation (\ref{prop:eau}) is written as $$K(t)= \dfrac{\tilde K}{[W](t)}. $$

%%%%%%%%%%%%%%%%%%%%%%%%%%%%%%%%%%%%%%%%%%%%%%%%%%%%%%%%%%%%%%%%%%%%%%%%%%%%%%%%%%%%%%%%%%%%
\subsection{Model 4}
This model is a mathematical simplification of the transport equation  
by means of homogenization theory. Homogenization theory is concerned  
with equations with rapidly oscillating coefficients and its aim is to  
provide an ``homogenized'' or ``averaged'' equation which is a  
limiting equation when the frequency of the oscillations tends to  
infinity. Advantages of  
homogenization theory are clear : on most occasions, it is simpler to  
use homogenized equations (for example to compute the  
solution) and, when the frequency of oscillations is over a given value, this  
approximation of the real equation by the homogenized one may be  
rather accurate as seen in the next section.

Homogenization problems for ODEs were studied by \cite{Pi}  
but it is worth pointing out that our particular case does not fall  
into the theory described in \cite{Pi}. Fortunately the specific  
structure of the transport equation allows us to do a complete   
analysis of the problem and even to compute explicitly the averaged  
equation.

More specifically, in the transport equation, pulses reach the bolus  
every 10 seconds approximately. Compared to the time scale of digestion
phenomena (the bolus stays in the small intestine for several hours \cite{wilfart,Beatrice Darcy et al.}), this represents a very high frequency and causes very rapid variations  
in the velocity of the bolus (see the velocity profile in Figure \ref{fig3.pdf}).

We can prove mathematically that the pulses can be
averaged out in an appropriate way and we can replace the rapidly  
varying velocity by a slowly
one. 

In the simplest case, by normalizing the pulses, we assume that their mean effect over a period is $e(\epsilon)$. Thus, over a time $t=N\epsilon$, their mean effect is $N e(\epsilon)=t e(\epsilon)/\epsilon$. On the other hand we assume
 $$\lim_{\epsilon \rightarrow 0}e(\epsilon)/\epsilon=\tau$$
 the mean effect over a time $t$ is therefore 
 $$\lim_{\epsilon \rightarrow 0}Ne(\epsilon)=te(\epsilon)/\epsilon=t\tau.$$

Inserting this equality in transport equation (\ref{equ:transport2}), the homogenized transport equation reads
\begin{eqnarray*}
&\dfrac{d^2x}{dt^2}(t) &= {\bar a}(t) \frac{c_0+c_1V(t)}{a+bx(t)}-\frac{\tilde K}{[W](t)} \frac{dx}{dt}(t)\\
&\dfrac{dx}{dt}(0) &= v_0 \label{prop:eau}, \, \, \, x(0)=0 
\end{eqnarray*}
where ${\bar a}$ is the averaged effect of the pulses. Its  
value is$$\displaystyle {\bar a}(t) :=
\tau (1-\frac{1}{c}\dfrac{dx}{dt}(t)).$$
%$\tau$ is the number of pulses per hour.
\section{Results }\label{result}
In the first part of this section, the graph of degradation of model 4, and the graph of transport of model 3 and 4 are developed.
The second part concerns the evaluation of the last model by comparing its outputs with experimental data. Only a limited number of outputs can be compared because of the lack of experimental data. However, the model is evaluated in relation to our objective which is developing a mathematical model that takes into account the physiology of the small intestine and process of digestion in it.
\subsection{Digestion}
The graph of digestion of model $4$ is shown in figure \ref{fig2.pdf}. We should at first initiate the bolus composition. These initial conditions vary following the different types of feedstuffs. 
We fixed the initial value of $A_{ns}$ as three times that of $A_{s}$. We dilute $A_{ns}$ by two times its volume water. Solubilized substrate $A_s$ and non-solubilized one $A_{ns}$ reach a dynamic equilibrium all along the small intestine, as explained in Section~\ref{digestion}. This balance is reached rapidly at the beginning of the small intestine due to the large difference in quantity between these two substrates. The result of this equilibrium is the increases of the value of $A_s$ and the  decreases of the values of $A_{ns}$, as seen in the graph of digestion. The inverse process might take place by lack of water. The absorption curve corresponds to the collected absorbable nutrients from $x=0$ to $x=x(t)$, where $x$ indicates the location of bolus in the small intestine. Obviously, the graph of the fourth model contains more details about different steps of digestion than the first two graphs thanks to the model structures. The digestion graph of model $4$ is similar that of model $3$. The only change in model $4 $ deals with the transport equation. 
\\[\intextsep]
\begin{minipage}{\linewidth}
  \centering%
 \includegraphics[width=140mm]{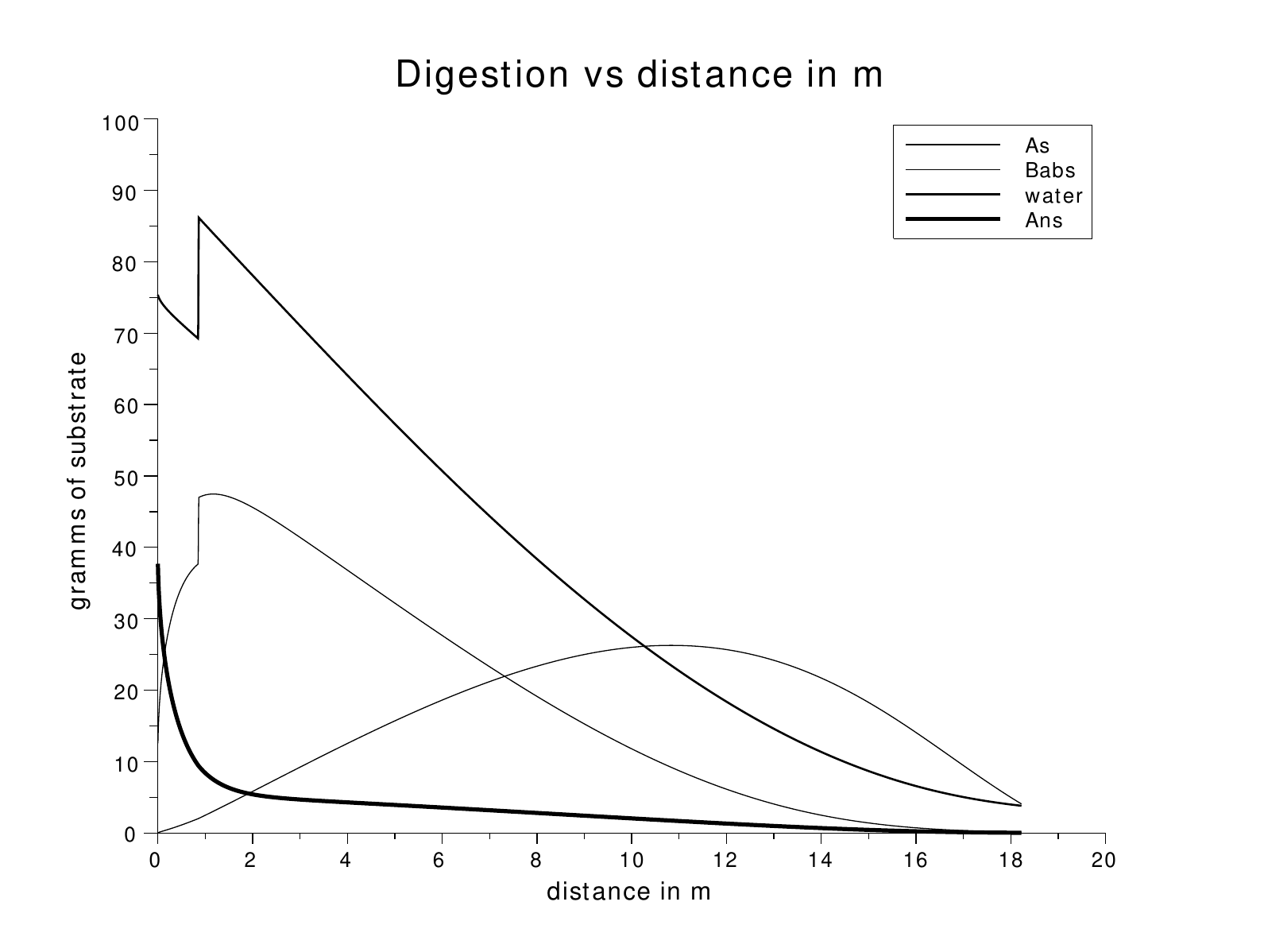}% 
\figcaption{Digestion through Model 4 }%
 \label{fig2.pdf}%
\end{minipage}
\\[\intextsep]

\subsection{Velocity}
Figure \ref{fig3.pdf} provides a numerical evidence of the homogenization phenomena. The graph of transport resulting from models $3$ and $4$ is shown in figure \ref{fig3.pdf}. The effect of pulses on the curve of velocity is obvious. However, using the homogenization theory in Model $4$, we obtain a smooth graph of velocity which replaced that of Model $3$.
\\[\intextsep]
\begin{minipage}{\linewidth}
  \centering%
  \includegraphics[width=140mm]{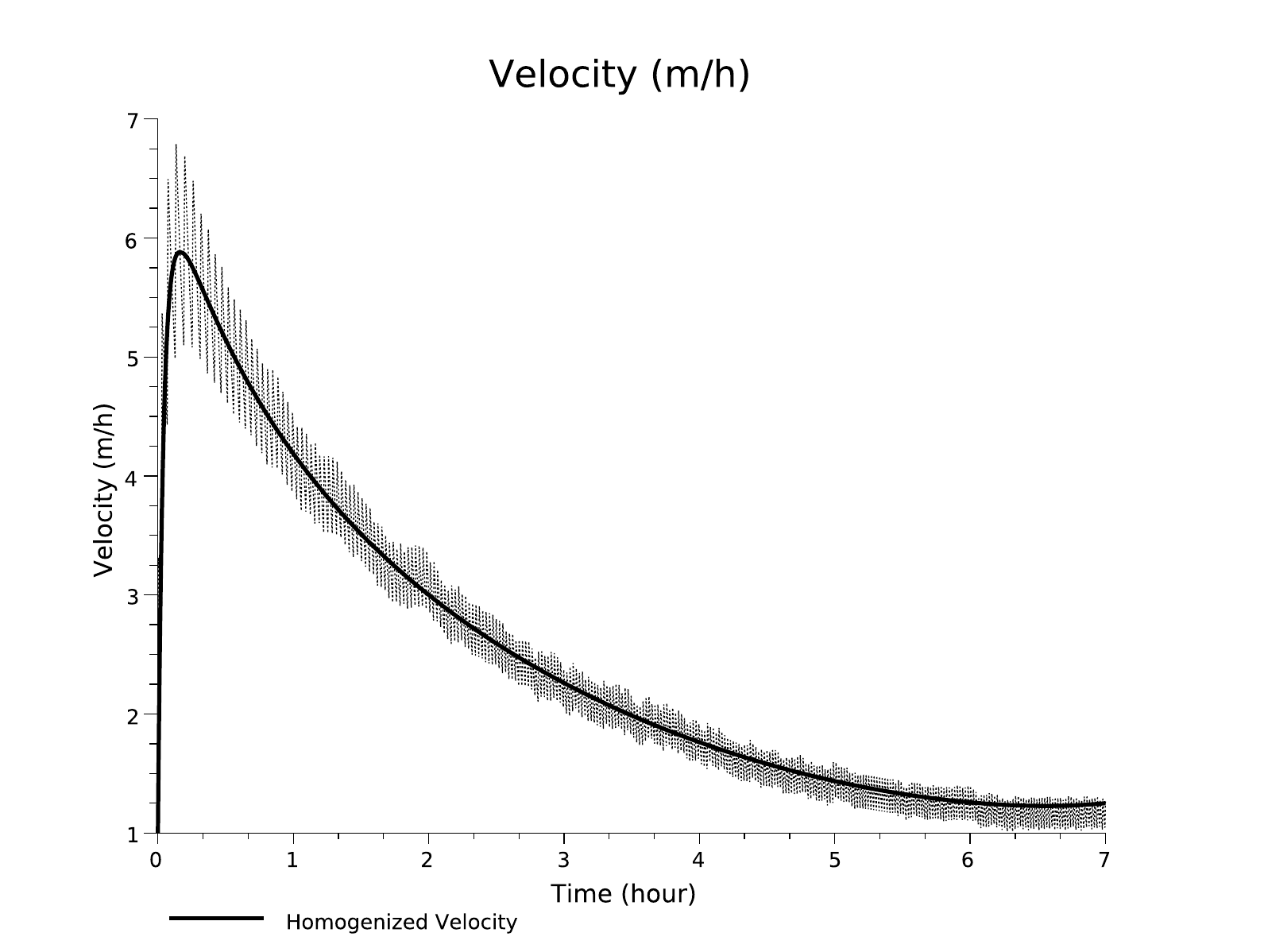}%
  \figcaption{Velocity of the bolus versus Time}%
  \label{fig3.pdf}%
\end{minipage}
\\[\intextsep]

\subsection{Model Evaluation}
 For a specific family of nutrients, here starch, digestion is calculated using Model 4 and is compared to data reported by \cite{Beatrice Darcy et al.}. To parameterize adequately the model, we adapt the enzyme activity of the last model to the activity of amylase in the small intestine. Amylase is the enzyme required for degradation of starch. The optimal activity of pancreatic amylase is in neutral pH \cite{digestion}.

Pancreatic secretions have no impact on the variation of $A_s$ since there is no starch from this source.

 The inputs of model are only $A_{ns}$ and $W$ which are Starch and Water. The outputs are the values of these substrates at the end of ileum. The data in the article of \cite{Beatrice Darcy et al.} are for purified protein free wheat starch, agreeing with our hypothesis for the composition of the bolus \ref{introduction}. The outputs concerns the collected data after at the end of ileum.  
\\[\intextsep]
\begin{minipage}{\linewidth}
  \centering%
\tabcaption{Digestion of Starch in Pigs by Modeling : Comparison between simulated and experimental data by \cite{Beatrice Darcy et al.}}
\label{tab:Commands}%
 \begin{tabular}{|c|c|c|}
\hline
 & Experimentation &  Modeling  \\
\hline
 \begin{tabular}{c}
\\
\hline
 wet digesta\\
\hline
dry matter  \\
\end{tabular}    &  \begin{tabular}{c|c}

Input(g)&Output(\%)\\
\hline
 $2571$& $8$  \\
\hline
 $688$  &  $0.50$(g)\\
                    \end{tabular}& \begin{tabular}{c|c}

Input(g)&Output(\%)\\
\hline
$113.10$ &  $5.33$  \\
\hline
$37.70$&   $0.04$ \\
                                        \end{tabular}\\
\hline
\end{tabular}
\end{minipage}
\\[\intextsep]

Regarding the data presented table \ref{tab:Commands}, percentages of dry matter and wet digesta collected at the end of the ileum are approximatively the same as the output of model $4$. The difference between inputs is due to the simulation calibration which takes into account only one bolus i.e. a fraction of the daily meal. However, differences between outputs are low in percentage enabling to conclude that the model can roughly simulate very simple situations.   
\subsection{Sensitivity Analysis}
Sensitivity analysis is performed to identify the key parameters affecting the digestion process. The chosen parameters are set at $5\%$ and $50\%$ of their original values.  
\\[\intextsep]
\begin{minipage}{\linewidth}
  \centering
%\label{tab:Commands}%
 \begin{tabular}{|c|c|}
\hline
Output & Parameters\\
\hline
 $A_s$  &  $C$, $C_{abs} $\\
\hline
$B_{abs}$ &   $C_{abs}$, $C_{iabs}$, $k_{abs}$  \\
\hline
$v$ & $a$,$b$, $c_0$, $c_1$, $K$\\
\hline
\end{tabular}
\end{minipage}
\\[\intextsep]

Studied digestion parameters are $C$ and $C_{abs} $ for degradation of $A_s$, and $C_{abs}$, $C_{iabs}$ and $k_{abs}$ for the absorption of $B_{abs}$. 

If $y$ is the output and $\theta$ the parameter, the relative variation of $y$ can be~expressed as follows
$$\frac{|y_{\theta}-y_{\theta +\Delta \theta}|}{y_{\theta}}.$$

\subsubsection{Influence on $A_s$ }

Both parameters $C$ and $C_{abs}$ are overestimated by $5$ and $50\%$. The figure \ref{fig:Assensitivity} shows the relative variation of $A_s$ at each moment. The relative variation of $A_s$ resulting from $5$ and $50\%$ values of $C$ is not meaningful. The parameter $C_{abs}$ has the largest effect on $A_s$ degradation. Observing the graph of relative variation of $A_s$, figure \ref{fig:Assensitivity}, we conclude that increasing the value of $C_{abs}$ increases the relative variation value with time. .
\\[\intextsep]
\begin{minipage}{\linewidth}
  %\centering
\begin{center}
 
\includegraphics[width=140mm]{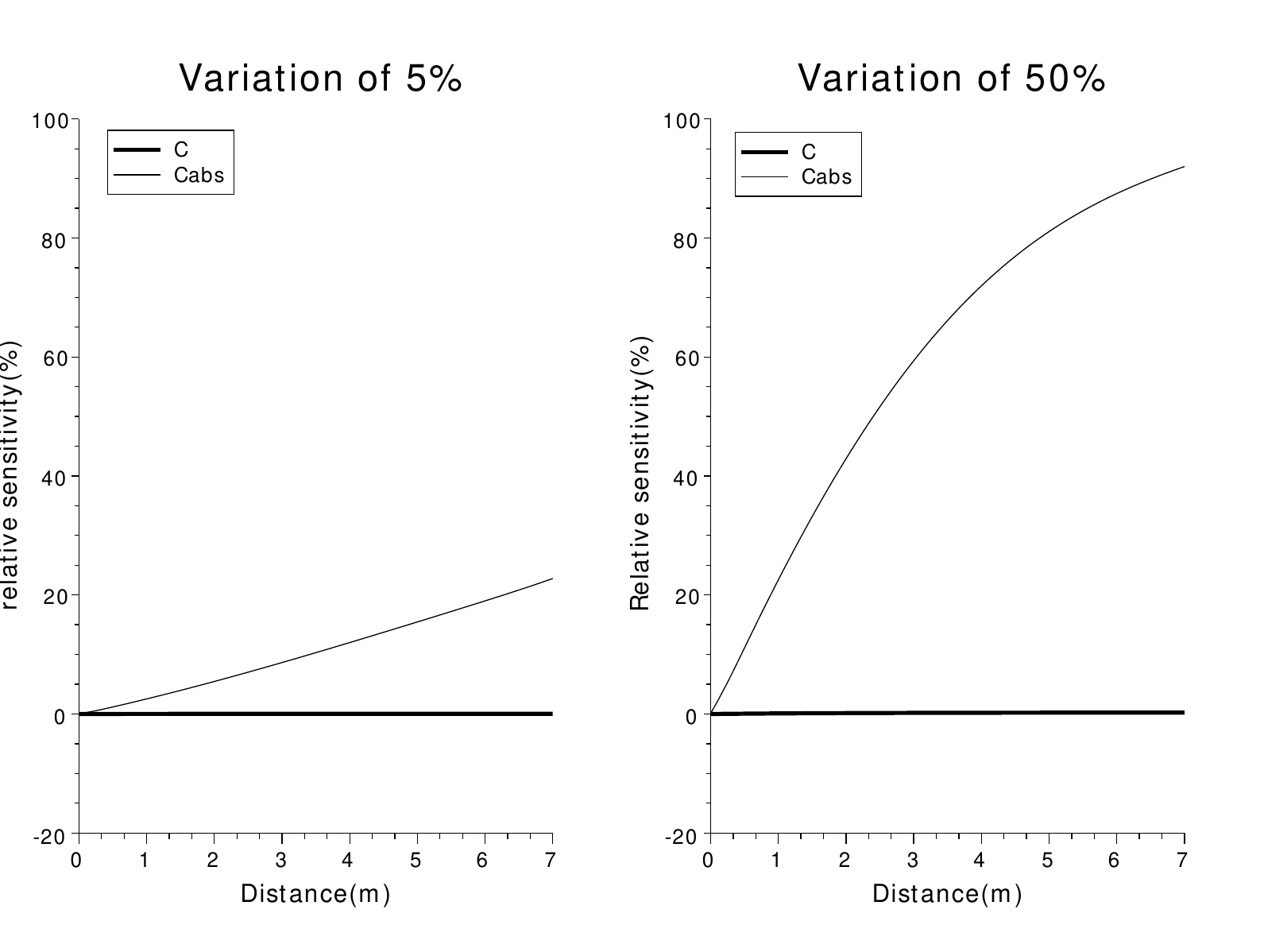}
  \figcaption{Relative variation of $A_s$ regarding to $C$, $C_{abs}$}%
\label{fig:Assensitivity}
\end{center}
\end{minipage}
\\[\intextsep]

\subsubsection{Influence on $B_{abs}$ }
The parameters $k_t$, $C_{abs}$ and $C_{iabs}$ are overestimated by $5$ and $50\%$. The figure \ref{fig:Babssensensitivity} shows the relative variation of $B_{abs}$ by time. The quantity of $B_{int}$ being very small in the model, the effect of changing the parameter $C_{iabs}$ is neglectible on the relative variation of $B_{abs}$ by time. The quantity $B_{abs}$ is very sensitive to the variation of parameter $C_{abs}$ firstly because of the high quantity of $A_s$, then its influence decreases because of decreasing quantity of $A_s$ over time. 
The quantity $B_{abs}$ is dependent on $kt$ because of the large impact of $kt$ on the nutrient absorption rate. 
\\[\intextsep]
\begin{minipage}{\linewidth}
  \centering
\includegraphics[width=140mm]{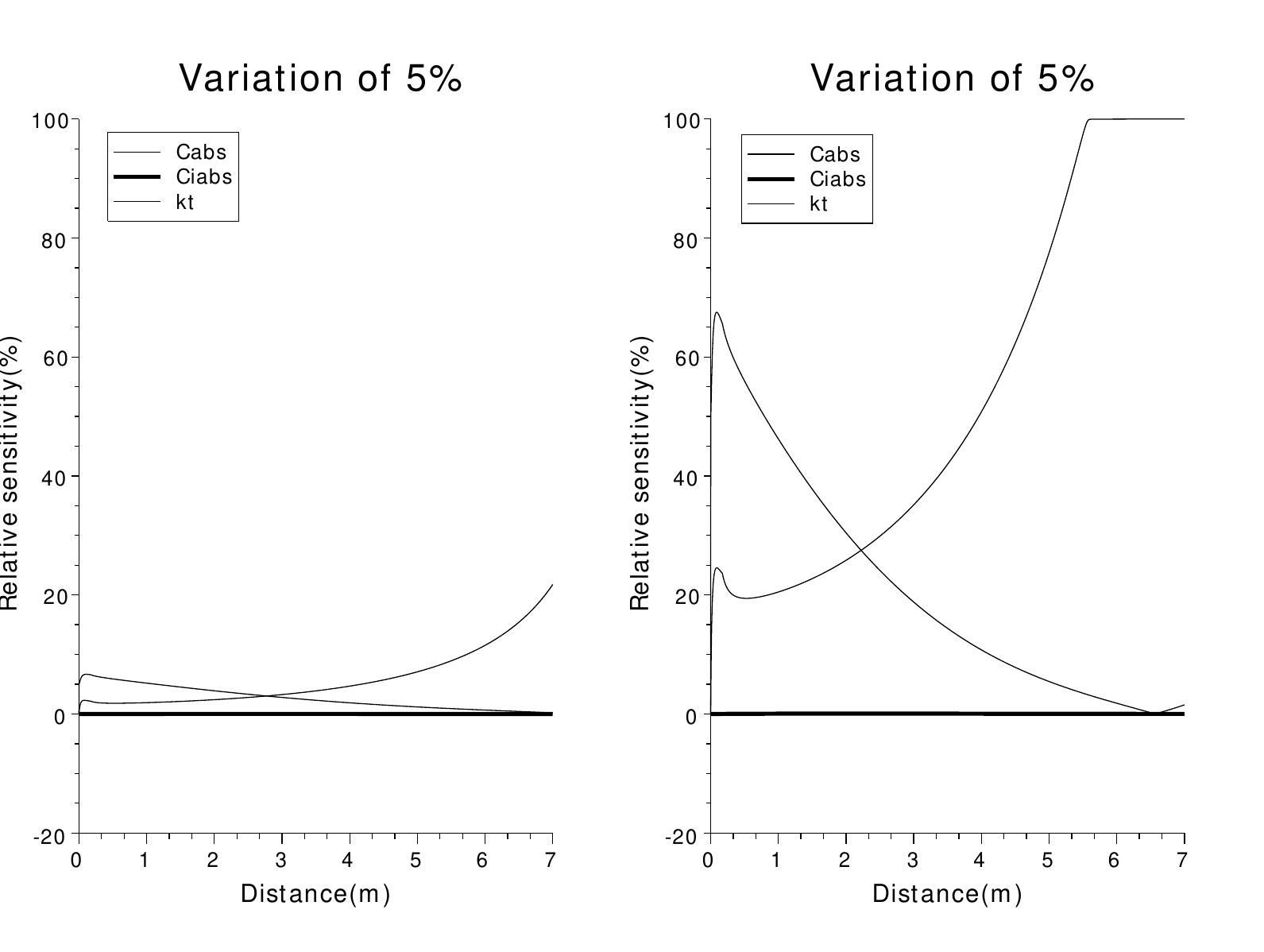}
 \figcaption{Relative variation of $B_{abs}$ regarding to $C_{abs}$, $C_{iabs}$ and $kt$}%
\label{fig:Babssensensitivity}
\end{minipage}
\\[\intextsep]

\section{Conclusion and Perspectives}\label{CandP}
This model is obtained from simplified biological assumptions and it can be used to illustrate generically the rate of degradation and absorption all along the small intestine. This is a global model of digestion of a bolus composed of one substrate and water. 
This section is devoted to a discussion on the current state of our modeling, our assumptions, difficulties and on the future development of the model.

The first assumption to be discussed is the ``cylinder'' one. As mentioned earlier, it was introduced for technical reasons. Our aim was to solve a partial differential equation with very different scales of times (pulses arising every $10$ seconds while the whole digestion process in the small intestine lasts for several hours \cite{wilfart, Beatrice Darcy et al.}) and with an highly variable domain with contrasted scales (few centimeters for the bolus compared to the $18$ meters of the small intestine). Solving this PDE seems unreasonable since it was leading to the usual diffusive phenomena and large errors. We also notice that these $18$ meters of the small intestine were rather empty and therefore we were spending a lot of time to compute functions which were very often $0$. 

From the transport standpoint, the ``cylinder'' assumption can be seen as a Lagrangian method, the ordinary differential equations on $x(\cdot)$ being (essentially) the {\em characteristic curves} of the transport equation. This is the first justification of this hypothesis, the second  being the direct observation of animals bolus which convinces us that it can be represented as a cylinder, even if its geometrical characteristics could be more complicated. However we have to work more on the evolution of the length of the cylinder.

A more Eulerian ``compartmental approach'' is studied simultaneously but we are still facing difficulties for modeling the transport and, specifically to reproduce some features of the ``cylinder'' model.

The transport equation seems to take into account rather closely the phenomena which are described by the experts . It will be difficult to validate the term $\displaystyle \frac{c_0 + c_1A}{a+bx}$ and to have a precise idea of the value of the different constants but such a modeling seems more appropriate than trying to use a complicated fluid mechanics approach whose laws may not be valid in this very confined domain. The same remarks hold for the effects of the water : it seems correct even if a relevant validation will be difficult.

For food digestion and absorption, we are only at a first stage of modeling. The absorption phenomena were not studied explicitly  leading to required further development with a focus on the assumed interactions between the animal physiological status and absorption. The spatial aspects (location of the absorption) were clearly neglected so far.

For digestion, the next step will be to mix different nutrients and adapt the enzyme breakdown to each of them. We have also to examine more closely the respective effects of the different categories of enzymes together with the role of the water. Moreover interactions between nutrients on the digestion processes should be questioned. 

As a conclusion of this first stage of modeling, consistent behaviors of the model  were reached. Moreover, the simplicity of the current model allows easy developments in any directions. Our next target will therefore be to iterate the model development according to the above proposed research areas.
\section{Acknowledgement}
The multidisciplinary collaboration on this research project between the INRA Center of Nouzilly and the Laboratoire de Math\'ematiques et Physique Th\'eorique was initiated within and supported by the CaSciModOT program (CAlcul SCIentifique et MOD\'elisation des universit\'es d'Orl\'eans~et de Tours) which is now a Cluster of the french Region Centre.
This collaboration also takes place in a CNRS-INRA PEPS program ``Compr\'ehension et Mod\'elisation du devenir de l'aliment dans le tube digestif``. This work is part of the PhD thesis of Masoomeh Taghipoor, financed by CNRS and INRA.

\end{document}